\documentclass[12pt]{amsart}
\usepackage{amsfonts,amssymb,amsmath}
\let\germ=\mathfrak

\let\Chi=\Xi
\newcommand{\Ind}{{\rm Ind}}
\newcommand{\Res}{{\rm Res}}

\let\eps=\varepsilon
 
\newtheorem{theorem}{\bf Theorem}

\newtheorem{lemma}[theorem]{\bf Lemma}
\newtheorem{definition}[theorem]{\bf Definition}
\newtheorem{proposition}[theorem]{\bf Proposition}
\numberwithin{equation}{section}
\title[Isolatedness of the minimal representation]
{Isolatedness of the minimal representation and minimal
decay of exceptional groups} 
\author{Hadi Salmasian}
\date{January 2 2005}
\address{Department of Mathematics, Yale University,
10 Hillhouse Avenue, New Haven, CT 
06520-8283.}
\email{hadi.salmasian@yale.edu}
\begin{document}
\begin{abstract}
Using a new definition of rank
for representations of semisimple groups,
sharp results are proved for the decay of matrix
coefficients of unitary representations 
of two types of non-split $p$-adic
simple algebraic groups of 
exceptional type.
These sharp bounds are achieved by minimal 
representations. It is also shown that in one of the
cases considered, the minimal representation is 
isolated in the unitary dual.
\end{abstract}

\maketitle
\section{introduction}

Let $\mathbf G$ be 
an absolutely simple
affine algebraic group defined over a local field $\mathbb F$
of characteristic zero. We assume $G$ is 
$\mathbb F$-isotropic.
Let $G$ be a topological finite 
central extension of 
$\mathbf G_\mathbb F$, the group of $\mathbb F$-points of $\mathbf G$.
Fix a maximal compact subgroup $K$ of $G$.
Let $\pi$ be a unitary representation of $G$. Throughout
this article we will
denote the Hilbert space on which $\pi$ acts by $\mathcal H_\pi$.
For any two vectors $v,w\in \mathcal H_\pi$, 
the matrix coefficient $f_{v,w}$ is the complex-valued 
function defined on $G$ as
$$
f_{v,w}(g)=(\pi(g)v,w)
$$
where $(\cdot,\cdot)$ is the inner product of $\mathcal H_\pi$.

\begin{definition}\label{stronglylp}
A unitary representation $\pi$ of $G$ is said to be
strongly $L^p$ if and only if for any vector $v$
in the (clearly dense) set of $K$-finite vectors
in $\mathcal H_\pi$, we have
$$f_{v,v}\in L^p(G,dg)$$
where $dg$ is the Haar measure on $G$.\\

We say that $\pi$ is strongly $L^{p+\eps}$
if and only if it is strongly $L^q$ for any
$q>p$.
\end{definition}

Cowling \cite{cowling} showed that $G$ has 
property $T$ of Kazhdan if and only if there exists 
a real number $p<\infty$, only dependent up on 
$G$, such that any non-trivial
irreducible unitary representation of $G$ is strongly $L^{p}$.
The infimum of all such numbers $p$ is called the 
{\it minimal decay} of $G$ \cite{LiZhu, Oh}. We denote 
this infimum by $p(G)$. The exact value of $p(G)$ 
is known for classical real Lie groups
\cite{Lidecay}. For
exceptional groups some estimates are obtained 
in \cite{Oh, LiZhu}. It is known that 
$p(G)=6$ when $G$ is of absolute type $\mathbf G_2$. 
Recently, Loke and Savin \cite{Savin}
used a notion of
$N$-rank developed by Weissman \cite{Weissman} for 
simply-laced groups to find $p(G)$ when $G$ is 
a split exceptional simply-laced group.\\

We use 
Tits' notation in \cite{Tits}  
to identify simple groups over $p$-adic fields.
The main point of this paper
is to prove the following statement:\\

\begin{theorem}
\label{mainthm}
$p({\!\!\ ^2\bf E_{6,4}^2})=8$, 
$p({\bf E_{7,4}^9})={26\over 3}$, 
and the minimal representation of
${\bf E_{7,4}^9}$
is isolated in its unitary dual 
endowed with the Fell topology.
\end{theorem}

\noindent{\bf Remarks.}\\
\noindent 1. We will see that the decay bounds
are
achieved by 
minimal representations of these groups.\\

\noindent 2. Although 
the notion of $N$-rank used by Savin is
defined in a way quite different from the one
introduced by the author,  
much of the methods used by
Savin and the author are similar, and based
on Howe's original ideas \cite{Howerank}.
The main 
result of \cite{Savin} follows from
the method of this paper
as well. On the other hand, the result of this paper
does not follow from the work of Loke and Savin
since there is no
parabolic with an abelian
unipotent radical when the relative root 
system is of type $\mathbf F_4$.  
However,  for the sake of brevity, 
we do not include the calculations
for split groups here.  
Besides, we address some issues about 
the minimal representations which are not addressed
in full generality in \cite{GanSavin}.
Therefore, this
paper is essentially a complement to the 
papers \cite{GanSavin} and \cite{Savin}.\\

\noindent {\bf Acknowledgement. }The 
author wishes to thank Roger Howe, 
the author's advisor,
and Gordan Savin for insightful 
communications.

\section{Rank of a unitary representation}

In this section we simply recall the 
definition of the notion of rank introduced in 
\cite{Salmasian} and the main result 
about purity of rank proved there.\\

Let $G$ be as in section 1. Throughout this section
we assume that 
$G$ satisfies properties $(H0)$ through $(H3)$ of
\cite[Section 3]{GanSavin}.
Take $G_1=G$ and 
let $P_1$ be the Heisenberg parabolic of
$G$, as defined in \cite{GanSavin},\cite{Weissman}
,\cite{Salmasian}. 
Let $P_1=L_1N_1$ be the Levi decomposition
of $P_1$. $N_1$ is a Heisenberg group unless 
$G_1$ is of type $\mathbf A_1$.\\

Note that in terms of relative root systems, 
the center of $N_1$ corresponds to 
the highest root $\beta$ of the relative 
root system
of $G$, and $N_1$
includes all unipotent subgroups
$U_\gamma$ correpsonding to relative roots $\gamma\in S$ 
where 

\begin{equation}
\label{defofS}
S=\{\gamma : (\gamma,\beta) > 0\}
\end{equation}

\noindent with $(\cdot,\cdot)$ being the Killing form.
There exists a unique simple
relative root $\alpha\in S$. \\

Suppose $L_1$ is isotropic. 
Then the relative root system of $L_1$
is either simple or of the form $\mathbf A_1\times \mathbf R$
or $\mathbf A_1\times\mathbf A_1\times\mathbf R$.
Let $G_2=[L_1,L_1]$  when $[L_1,L_1]$ is simple, and 
otherwise let $G_2$ be equal to the simple
factor of $[L_1,L_1]$ which corresponds
to the system $\mathbf R$,
and let $P_2=L_2N_2$ be the Heisenberg 
parabolic of $G_2$. Again if $L_2$ is isotropic we 
can consider its simple factor $G_3$ and
the Heisenberg parabolic $P_3=L_3N_3$ 
of $G_3$, and so on.
The unipotent subgroup

\begin{equation}\label{OKPgroupN}
N_G=N_1\cdot N_2\cdot N_3\cdots
\end{equation}

\noindent of $G$ is in fact the unipotent radical 
of a parabolic subgroup
of $G$. $N_G$ is a tower of extensions by
Heisenberg groups; however, the last group in the 
sequence can be an abelian group.\\

One can also construct a family of 
representations of $G$ in a similar fashion.
Any irreducible  
unitary representation
of a Heisenberg group over a local field belongs
to one of the following classes. The first class 
consists of characters, i.e. representations which
act trivially on the center of the group. 
The second class
consists of those representations which have a 
non-trivial central character. By the Stone-von Neumann
theorem, the elements of
this class are uniquely determined by their central 
characters. See \cite{Howeheisenberg} or 
\cite{Taylor} for more details.\\

Now let $\germ H$ be a Heisenberg group, and 
let $\chi$ be an arbitrary non-trivial
character of the center of $\germ H$.
Denote the unitary representation with this central
character by $\rho_\chi$. \\

\begin{definition}
\label{definition_of_rankable}
(\cite[Section 4]{Salmasian})\ \ 
Let $N_G$ be as in (\ref{OKPgroupN})
with height $r=r_{N_G}$; i.e. 
$$N_G=N_1\cdot N_2\cdots N_r.$$
Let $\rho_{\chi_i}$ be the extension of 
a representation of $N_i$ with central character
$\chi_i$ to $N_G$
\footnote{$\rho_{\chi_i}$ is first extended to $N_{i+1}\cdots N_r$ using
the Weil representation \cite{Weil} and then trivially on the normal
subgroup $N_1\cdots N_{i-1}$.}. Any representation
of the form 
$\rho_{\chi_1}\otimes\rho_{\chi_2}\otimes\cdots
\otimes\rho_{\chi_k}$ is called a rankable representation
of $N_G$ of rank $k$.\\

%\noindent The height $r$ of $N_G$
%is denoted by $r_{N_G}$.
\end{definition}
The following proposition is a special, more
polished form of the main result in 
\cite[Section 5.2]{Salmasian}.

\begin{proposition}
Let $G$, $N_G$ be as in Definition \ref{definition_of_rankable}. 
Suppose $r_{N_G}>2$.
Let $\pi$ be a unitary representation of $G$ without a $G$-fixed 
vector.
Then the spectral measure 
for the direct integral decomposition of 
$\pi_{|N_G}$ is supported on
rankable unitary representations. 
Moreover, if $\pi$ is irreducible,
then exactly one of the following cases can occur:

\begin{enumerate}
\item[i.] The spectral measure is supported on rankable
representations of rank one.
\item[ii.] The spectral measure is supported on rankable 
representations of rank strictly larger than one.
\end{enumerate}

\end{proposition}

\noindent Therefore, the following definition of rank, 
given in \cite[Definition 5.2.3]{Salmasian}, 
distinguishes the representations of rank one
from other representations.

\begin{definition}\label{rank_for_representations}
Let $G$ be as in section 1,
$N_G$ be as in (\ref{OKPgroupN}), and 
$\pi$ be an irreducible unitary representation of
$G$. $\pi$ is said to have rank $k$ if and only if 
the spectral measure in the direct integral
decomposition of $\pi_{|N_G}$ is supported on 
rankable representations
of $G$ of rank $k$.
\end{definition}

\section{Matrix coefficients}

Let $G$ be as in section 1 with maximal
compact subgroup $K$. 
Let $B$ be the minimal parabolic of $G$
and $A$ the split torus of $G$.
Since $G=BK$, for any $g\in G$ there exists
an element $b=b(g)\in B$ 
(unique modulo $B\cap K$)
such that $b^{-1}g\in K$.
Let $\delta_G$ be the modular function of $B$.
The Harish-Chandra function $\Chi$, which is the 
spherical matrix coefficient of 
$\mathrm{Ind}_B^G(1)$, is equal to 
\begin{equation}\label{shaifunction}
\Chi(g)=\int_K\delta_G^{1\over 2}(b(kg))dk. 
\end{equation}
Note that $\delta_G(a)$ assumes values
in $\mathbb R^+$, so it is identically equal
to 1 on compact 
subgroups of $B$; i.e. $\delta_G(b(g))$ is 
well defined.\\
\begin{definition}
Let $\pi$ be a unitary representation of $G$.
Let $K$ be the maximal compact subgroup of $G$,
and let $\Psi$ be a positive integer-valued function
whose domain is equal to the 
set $\hat K$ of the equivalence classes of 
finite-dimensional irreducible representations of
$K$.\\

\noindent We say $\pi$ is $(\Chi^{1\over k},\Psi)$-bounded
if and only if for any $K$-finite vectors 
$v, w\in\mathcal H_\pi$
which belong to $K$-isotypic spaces of $\nu_v, \nu_w\in \hat K$,
we have 
$$
|f_{v,w}(g)|\leq \Chi^{1\over k}(g)\Psi(\nu_v)\Psi(\nu_w).
$$ 
\end{definition}
Let $t>1$. It is known (e.g. \cite[Prop. 5.2.8]{Knapp}
and \cite[Lemma 4.1.1]{Silberger})
that if 
$$
A_t=
\{a\in A : |\alpha(a)|\geq t
\ \textrm{for any positive root}\ \alpha\}
$$
then for any $a\in A_t$ we have
\begin{equation}\label{deltainequality}
C\delta_G(a)^{1\over 2}\leq 
\Chi(a)\leq C_t\delta_G(a)^{1\over 2}.
\end{equation}
\begin{proposition}
\label{p6}
Let $G$ be as in section 1 and $K$ be 
the maximal compact of $G$. Let $\Chi$ be 
defined as in (\ref{shaifunction}), 
and $\Psi(\nu)$ be
an arbitrary function from the $K$-types to 
$\mathbb Z^+$.\\
\begin{enumerate}
\item[1.] The subset of $(\Chi^{1\over k},\Psi)$-bounded
representations of $G$ is closed in the unitary dual of 
$G$. Any such representation $\pi$ is strongly 
$L^{2k+\eps}$. 
%In fact, all $K$-finite matrix coefficients
%of $\pi$ are in $L^{2k+\eps}(G,dg)$ for any $\eps>0$.
\item[2.] Let $\pi$ be a unitary representation of 
$G$ which is  
strongly $L^{p+\eps}$. Let $k$ be a positive
integer such that $p\leq 2k$.
Then 
$\pi$ is $(\Chi^{1\over k},\mathrm{dim}(\nu))$-bounded. 
\end{enumerate}
\end{proposition}
\begin{proof}
These are restatements of Howe's results in
\cite{Howerank} in general form.
Use (\ref{deltainequality})
for the first part. For the second part
see \cite[Corollary 7.2]{Howerank}.
\end{proof}
Let $H_1\times\cdots\times H_m\subset G$ be a product
of semisimple groups embedded 
inside $G$, such that
if $A_i$ denotes the maximal split torus of $H_i$,
then $A_1\times\cdots\times A_m$ is equal to the maximal
split torus $A$ of $G$.
Choose positive Weyl chambers $A^+,A_1^+,...,A_m^+$
such that 
$A_1^+\times\cdots\times A_m^+\subseteq A^+$.
Let $\delta_{H_i}$ be the modular function
of the minimal parabolic subgroup of $H_i$. Obviously
we can assume $\delta_{H_i}$ is also defined
on all of $A$ so that it is the modular function of
the minimal parabolic of the reductive group 
$A\cdot H_i$.

\begin{proposition}
\label{asledecay}

Assume the setting introduced above.
Let $\pi$ be a unitary representation of 
$G$ such that $\pi_{|H_i}$ is strongly
$L^{p_i+\eps}$. 
Let $k_i$'s be positive integers such that
$2k_i\geq p_i$. 
Suppose 
$$
\textrm{for any \ } a\in A,
\qquad\prod_{i=1}^m\delta_{H_i}(a)^{1\over 2k_i}\geq 
\delta_G(a)^{1\over p}.
$$
Then $\pi$ is strongly $L^{p+\eps}$.
\end{proposition}
\begin{proof}
This follows from \cite[Theorem 3.1]{Lidecay} or 
\cite[Prop. 2.7]{Oh}.
\end{proof}

\section{Minimal representations of archimedean groups}

In this section and the next one 
we address the classification of rank one
and minimal representations. The classification of 
such representations in the classical case
is already known by the work of Howe, J.-S. Li and Vogan. 
Here we study exceptional groups. Therefore, in this
section we assume 
$G$ is as in section 2, but also of exceptional type.\\

The traditional
definition of minimal representations
is different for the archimedean and non-archimedean
case. A notion 
of a weakly minimal representation is introduced in
\cite[Definition 3.6]{GanSavin} in the non-archimedean case.
In fact, that definition simply means that the representation
is of rank one in the sense of Definition 
\ref{rank_for_representations}
above. 
Therefore it is quite natural
to generalize it to the archimedean
case as well. In fact, as seen below, 
this new definition of
the
minimal representation agrees
with the traditional one which is 
in terms of the Joseph
ideal.\\

Let $G$ be as in section 2, and let $\pi$ be a unitary 
representation of $G$. Suppose $P$ is the Heisenberg 
parabolic of $G$. By Mackey theory
\cite{Mackey}, 
$\pi_{|P}$ is 
a direct integral of representations of the form

\begin{equation}
\label{mackeyPP}
\Ind_{[P,P]}^P(\nu_\chi\otimes\rho_\chi)
\end{equation}

\noindent where $\nu_\chi$ is a representation of $[P,P]$ which
factors through $[P,P]/N$, and $\rho_\chi$ is 
an extension of the representation of $N$ with 
central character $\chi$ to $[P,P]$. 
Note that $\chi$ is supposed to be a nontrivial 
unitary character.
Since the action of 
the split torus of $G$ on the center of $N$ has only
finitely many orbits, $\pi_{|P}$ will be a 
finite direct sum of representations of the 
form given above.\\

We now recall the definition of the minimal 
representation in the archimedean case. For
simplicity let us assume 
$\mathbb F=\mathbb R$. The 
complex case can be treated similarly. 
Let
$\germ g_\circ$ be the Lie algebra of $G$, with
complexification 
$\germ g=\germ g_\circ\otimes\mathbb C$.
$\germ g$ has a complex 
Heisenberg parabolic subalgebra
$\germ p$ with 
Levi decomposition 

\begin{equation}
\label{LeviofHeisenberg}
\germ p=\germ l\oplus\germ h.
\end{equation}

\noindent Note 
that $\germ h$ will be the complexification
of the Lie algebra of $N$ where $P=LN$.
Let $\germ m=[\germ l,\germ l]$.
Assume $\germ g$ is not of type $\mathbf A_l$, 
and let $\mathcal J_\circ$ be the Joseph ideal of
$U(\germ g)$. For the definition of the 
Joseph ideal, see \cite{Joseph} or \cite{GanSavin-unique}. \\

Let $\mathcal U(\germ h)$ be the universal
enveloping algebra of $\germ h$.
Define the map 
$$
\Theta: \germ m\mapsto\mathcal U(\germ h)
$$
as follows
\footnote{This neat description is 
borrowed from \cite{GanSavin}.}. Let
$$
\germ h=W\oplus\germ z
$$
where $W$ is a symplectic space and $\germ z$ is the center
of $\germ h$. Denote the symplectic bilinear form
of $W$ by $\omega(\cdot,\cdot)$.
We can identify $\germ{sp}_n$
with the symplectic Lie algebra $\germ{sp}(W)$.
There is a canonical map
$$
\germ{sp}(W)\mapsto (W^*\otimes W^*)^{S_2}
$$
where $S_2$ is the symmetric group on two elements
( i.e.
$S_2=\{\pm1\}$), which is defined by sending
an element $A\in\germ{sp}(W)$ to the bilinear form
$$
\omega_A(x,y)=\omega(Ax,y).
$$
Using $\omega$, we can identify $W^*$ 
with $W$ by sending any $y$ to the linear
form $y\mapsto \omega(x,y)$. Thus we have a 
map 
$$
\germ{sp}(W)\mapsto (W\otimes W)^{S_2}.
$$
Since there is a natural map
$(W\otimes W)^{S_2}\mapsto \mathcal U(\germ h)$, we obtain
a composition map $A\mapsto n_A$ from $\germ{sp}_n$ to $\mathcal U(\germ h)$.
Now set 
$$
\Theta(X)={1\over 2}n_X.
$$
We have
$$
[\Theta(X),\Theta(Y)]=Z\Theta([X,Y]).
$$

By \cite[Prop. 3.1.1]{Salmasian}, 
we can consider 
$\germ m$ as a Lie 
subalgebra of 
a complex symplectic 
Lie algebra $\germ{sp}_n$ which acts in 
the usual way on the
Heisenberg algebra $\germ h$  
(see \cite[Section 1]{Howerank}). Thus we 
can define $\Theta(X)$
for any $X\in\germ{m}$.

\begin{lemma}
\label{KOSTANTGARFINKLE}
Let $\mathcal J$ be a primitive ideal
of infinite codimension
in $\mathcal U(\germ g)$. Let $\germ l$ be
as in (\ref{LeviofHeisenberg}) and
$\germ m=[\germ l,\germ l]$.
Suppose
$\mathcal J$ contains
$ZX-\Theta(X)$ for all $X\in\germ m$.
Then $\mathcal J$ is the Joseph ideal.
\end{lemma}
\begin{proof}
This follows immediately from 
\cite[Proposition 4.3]{GanSavin} 
and the uniqueness of the Joseph ideal
\cite{GanSavin-unique}.
\end{proof}

\begin{definition}
\label{classicalminimal}
Let $G$ be as above (i.e. $F=\mathbb R$).
An irreducible
unitary representation $\pi$ of $G$ is called
minimal iff the annihilator of the Harish-Chandra
module associated to $\pi$ is the Joseph ideal.
\end{definition}

\begin{proposition}
Let $\pi$ be an irreducible 
unitary representation of $G$.
$\pi$ is a minimal representation if and only if 
$\pi$ is rank one (in the sense of Definition
\ref{rank_for_representations}).
\end{proposition}

\begin{proof}
Let $\pi$ be a representation of $G$ 
of rank one. We will show that $\pi$ is
a minimal representation.
Intuitively,
the proof is as follows. 
Let $d\pi$ denote the 
infinitesimal $\germ g$-action. 
By Mackey theory, the restriction of $\pi$
to $G$ is a finite direct sum of 
representations given in (\ref{mackeyPP}).
But the representations $\nu_\chi$ will
be of rank zero, i.e. they are trivial
$[P,P]$ modules. Consequently, 
for some $m$, we have

\begin{equation}
\label{wplusminus}
\pi_{|P}=\bigoplus_{i=1}^m
n_i\Ind_{[P,P]}^{P}\rho_{\chi_i}
\end{equation}

\noindent where $n_i\in\{0,1,2,...,\infty\}$ for each $i$.
Note that in fact $m=1$ except for the Hermitian case, where
$m=2$. Let

\begin{equation}
\varpi_i= 
\Ind_{[P,P]}^{P}\rho_{\chi_i}.
\end{equation}

\noindent By  
\cite[Lemma 4.3.2]{Salmasian} 
or \cite[Theorem 5.19]{grahamvogan}
it follows that for any nontrivial unitary character
$\psi$,

\begin{equation}
\label{rhokills}
d\rho_\psi(ZX-\Theta(X))=0\quad\textrm{for any}
\quad X\in \germ m
\end{equation}

\noindent and therefore a similar 
equality should be true for 
each $d\varpi_i$ and also for $d\pi$; i.e.

\begin{equation}
\label{imprecisekill}
ZX-\Theta(X)\in 
\mathrm{Ann}_{\mathcal U(\germ g)}(\pi).
\end{equation}
The result now follows
from Lemma \ref{KOSTANTGARFINKLE}.\\
    
We now have to give a rigorous proof
of (\ref{imprecisekill}).
Let    
$\mathcal H_\pi$ be the Hilbert space of the 
representation $\pi$, and $\mathcal H_\pi^\infty$ be
its subset of $G$-smooth vectors. \\

\noindent Let $\varpi\in\{\varpi_1,...,\varpi_m\}$. We have
the following lemma.
\begin{lemma}
\label{KILLP}
Let $\mathcal H_\varpi$ denote the 
Hilbert space of the representation
$\varpi$, and let $\mathcal H_\varpi^\infty$
be the set of $P$-smooth 
vectors in 
$\mathcal H_\varpi$. Then 
$$
d\varpi(ZX-\Theta(X))v=0\quad
\textrm{for any}\ X\in\germ m,v\in
\mathcal H_\varpi^\infty.
$$
\end{lemma}
\begin{proof}
Recall that $\mathcal H_\varpi^\infty$ is  
a space of smooth left-quasi-$[P,P]$-invariant
functions on $P$ where $P$ acts on it by the 
right regular representation.
Let $X\in\germ m$. Take any 
$f\in\mathcal H_\varpi^\infty$. 
Let $e$ be 
the identity element of $P$, 
and denote the right action of $p\in P$ 
on
$\mathcal H_\varpi$ by $R_p$. Thus
$\varpi(p)$ is equal to $R_p$.

%a normalization
%of $R_p$ (because $P$ is not unimodular).\\
	
Let $M_X=ZX-\Theta(X)$.
For any $p\in P$
we have
\begin{align}
\begin{split}
\notag
	\left(d\varpi(M_X)f\right)(p)
	=(d\varpi(Ad(p)(M_X))(R_pf))(e)\\
	=d\rho_{\psi}(Ad(p)(M_X))
	(f(p))
	=d\rho_{p^{-1}\cdot\psi}
	(M_X)(f(p)).
	\end{split}
	\end{align}
    Here $p\cdot\psi(x)=\psi(p^{-1}xp)$  
    as usual. An application of
   (\ref{rhokills}) completes the
    proof.
\end{proof}
By Lemma \ref{KILLP}, 
$d\varpi(ZX-\Theta(X))v=0$ for any 
$v\in \mathcal H_\varpi^\infty$. Since $\pi_{|P}$
is a (finite or infinite) direct sum of representations
unitarily equivalent to $\varpi$, the algebraic sum of 
the spaces of $P$-smooth vectors will be a 
norm-dense $P$-invariant subspace of 
the space $\overline{\mathcal H}_\pi^\infty$
of $P$-smooth vectors in $\mathcal H_\pi$.
However, by a result of 
Poulsen \cite{poulsen}, this space is dense in the 
Fr\'{e}chet topology of 
$\overline{\mathcal H}_\pi^\infty$,
and therefore
$d\pi(ZX-\Theta(X))v=0$ for any 
$v\in\overline{\mathcal H}_\pi^\infty$. 
But $\mathcal H_\pi^\infty\subseteq 
\overline{\mathcal H}_\pi^\infty$, which 
proves that $\pi$ is a minimal representation
in the sense of Definition \ref{classicalminimal}.\\

Now we address the converse.
Let $\pi$ be a minimal representation
in the sense of Definition \ref{classicalminimal}. 
We think of $\pi_{|P}$ as
a representation of $\tilde P^\circ$, 
the universal cover of the connected
component of identity of $P$.
Again by Mackey
theory, we can express $\pi_{|\tilde P^\circ}$ 
as
$$
\pi_{|\tilde P^\circ}=\bigoplus_{i=1}^m 
\nu_{\chi_i}\otimes
\Ind_{[\tilde P^\circ,\tilde P^\circ]}^{\tilde P^\circ}\rho_{\chi_i}
$$
where now $\nu_{\chi_i}$'s
will be representations of $\tilde L^\circ$.
If $\pi$ is minimal, then Lemma \ref{KILLP}
implies that the actions of $d\nu_{\chi_i}$'s
have to be zero; i.e. $\nu_{\chi_i}$'s 
have to be trivial modules.
This proves that $\pi$ is a rank one
representation.\\

\end{proof}

\section{Classification of rank one representations}

Let $G$ be as in section 2 and also assume 
$\mathbf G$ is
of exceptional type.
Let $P=LN$ be the Heisenberg 
parabolic of $G$. Let $\pi$ be an irreducible  
representation of $G$ of rank one. Then we have a
decomposition as in (\ref{wplusminus}) 
with $\nu_{\chi_i}$'s being
trivial modules. Let $s_\alpha$ be the reflection
inside the relative Weyl group of $G$ which 
corresponds to the simple root $\alpha\in S$ 
(see (\ref{defofS}) ). By Mackey's subgroup 
theorem \cite{Mackey} one
can see that

\begin{lemma}
\label{BsBs}
If $\varpi_i$ and $\varpi_j$ correspond to
different orbits of the split torus on the 
center of $N$, then
$\Res_{P\cap s_\alpha Ps_\alpha^{-1}}^P\varpi_i$ 
and
$\Res_{P\cap s_\alpha Ps_\alpha^{-1}}^P\varpi_j$
are non-isomorphic  
irreducible representations of 
the group $P\cap s_\alpha P s_\alpha^{-1}$.
\end{lemma}

Let $\pi$ be an irreducible representation
of $G$ of rank one.
From the 
irreducibility claimed in 
Lemma \ref{BsBs} it follows that if $T$
is a bounded operator from the space of $\pi$ to itself and
it commutes with all $\pi(g)$ for 
$g\in P\cap s_\alpha P s_\alpha^{-1}$, then 
it commutes with all $\pi(g)$ for $g\in P$ as well.
For a similar reason, $T$ should commute with all 
$\pi(g)$ for $g\in s_\alpha P s_\alpha^{-1}$ 
as well. However, $P$ and $s_\alpha P s_\alpha^{-1}$
generate $G$; therefore $T$ should be an intertwining operator
of the irreducible representation $\pi$, which means that
$T$ is a scalar. Consequently, we have proved that
the only intertwining operators of 
$\pi_{|P\cap s_\alpha P s_\alpha^{-1}}$ are scalars.
This means that $\pi_{|P}$ should be irreducible; i.e.
$\pi_{|P}=\varpi$ for some $\varpi=\varpi_i$. Moreover,
$\pi(s_\alpha)$ will be the intertwining operator between
two irreducible 
representations of $P\cap s_\alpha Ps_\alpha^{-1}$, and therefore
$\pi(s_\alpha)$ will be determined uniquely. Therefore once
we know the restriction of $\pi$ to $P$, there is at most
one way to extend it to a representation of $G$.\\

\begin{proposition}\label{UNIQUENESSOF}
Let $G$ be as above. Suppose  
the absolute root system of
$G$ is of exceptional type and
$G$ has split rank at least 3.
Then there is a unique
representation of rank one of $G$, 
unless $G$ is a Hermitian form of
$E_7$ (i.e. $\mathbb F=\mathbb R$), 
in which case 
there are exactly two such representations.
\end{proposition}

\begin{proof}
From the discussion above, and the fact that 
representations $\varpi_i$ correspond to distinct orbits of
the center of $N$, 
the proposition 
follows from the number of orbits of the split
torus on the center of $N$. Consider the smallest subspace
of the Lie algebra of $G$ which contains 
the one-dimensional space of the 
highest root and is invariant under the action
of the rank-one subgroup corresponding to $\alpha$. 
It is either the two-dimensional 
representation of $SL_2(\mathbb F)$ or a 10-dimensional representation
of $SO(9,1)$ (in the 
real Hermitian case). 
In the former case, the split torus 
of $SL_2(\mathbb F)$ has one orbit, 
so there can be at most one minimal representation
in each case. 
The representations
are constructed in \cite{torassoduke}.
In the latter case,
there are two orbits. The two minimal representations
are constructed in \cite{Sahi}.

\end{proof}

\section{Decay of minimal representations}

As we mentioned in the introduction, the 
method of this paper can be used to show that the
minimal representation is isolated for all 
simply-laced exceptional groups as well as two
non-split exceptional $p$-adic groups. However, 
the details for the case of split groups is not included here
since the author noticed that it 
was done simultaneously 
in \cite{Savin}. \\

The non-split $p$-adic groups under consideration
will be the $\mathbb F$-points of simply connected
algebraic groups of absolute 
types $\mathbf E_6$ and $\mathbf E_7$
which
have Tits index $\!\!\ ^2\mathbf E_{6,4}^2$ and $\mathbf E_{7,4}^9$.
They are of relative type $\mathbf F_4$, and the dimension of  
short root spaces are $2$ and $4$ respectively. 
By Proposition \ref{UNIQUENESSOF} the minimal
representation is unique in each case, and all other
nontrivial irreducible  
unitary representations are of rank two or larger. 
The minimal representations have Iwahori-fixed vectors and the 
exponents of these representations are calculated in 
\cite[Section 8]{GanSavin}. By a standard result relating the
matrix coefficients to the exponents 
(e.g. see \cite[Section 4]{Casselman}), one can calculate the 
precise $L^p$ decay of the minimal representations.
Let the restricted Dynkin diagram of $G$ be labelled 
by $\alpha_1,...,\alpha_4$ 
such that the highest root $\beta$ is

$$
\beta=2\alpha_1+3\alpha_2+4\alpha_3+2\alpha_4.
$$
 
\noindent Note that this is different from the labelling
chosen
by Lusztig and also used in \cite{GanSavin}. 
Let us choose a minimal parabolic $B$ of $G$, 
and therefore identify the positive
Weyl chamber of the split torus $A$.\\

Let $\pi_{\min}$ denote the (unique) 
minimal representation of $G$, as constructed in
\cite{GanSavin}. We use the calculation of exponents
of $\pi_{\min}$ given in \cite[Section 8]{GanSavin}.
See \cite[Section 6]{GanSavin} for the definition 
of and 
basic results on exponents.
The following statement follows from
standard facts (see \cite[Section 4]{Casselman}
or \cite{HGL}).
\begin{center}
\begin{itemize}
\item[$\Diamond$]
$\pi_{\min}$ is strongly $L^{p+\eps}$ if and only if
for every exponent $\mu$ of $\pi_{\min}$,
$\delta_G^{{1\over 2}-{1\over p}}(a)\mu(a)$ is bounded on the 
negative Weyl chamber.\\
\vspace{0.2cm}
\end{itemize}
\end{center}

\noindent 
We treat each case separately below.\\

\begin{enumerate}
\item[$\bullet$] $\ ^2\mathbf E_{6,4}^2$:\\

\noindent By \cite[Section 8.4]{GanSavin},
the logarithms of absolute values
of exponents of $\pi_{\min}$ are

\begin{equation}
\label{exponentsofe6}
\left\{
\begin{array}{l}
\log|\overline{\chi_4}|=
-8\alpha_1-15\alpha_2-22\alpha_3-12\alpha_4\\
\log|\chi_4|=
-8\alpha_1-15\alpha_2-22\alpha_3-12\alpha_4\\
\log|\chi_2|=
-7\alpha_1-15\alpha_2-22\alpha_3-12\alpha_4
\end{array}
\right.\end{equation}

\noindent Moreover, 

$$
\log \delta_G=22\alpha_1+42\alpha_2+60\alpha_3+32\alpha_4.
$$

$\pi_{\min}$ is strongly $L^{p+\eps}$ if and only if
for any $\chi$ in \ref{exponentsofe6},
the coefficients of $\alpha_i$'s 
in 
$$
\chi+ ({1\over 2}-{1\over p})\log \delta_G
$$
are nonnegative. This is equivalent to $p\geq 8$.\\

\item[$\bullet$] $\mathbf E_{7,4}^9$ :\\

\noindent The logarithms of absolute values
of exponents are equal to

\begin{equation}
\left\{
\begin{array}{l}
\log|\chi_4|=-13\alpha_1-24\alpha_2-36\alpha_3-20\alpha_4\\
\log|\chi_2|=-11\alpha_1-24\alpha_2-36\alpha_3-20\alpha_4
\end{array}
\right.
\end{equation}

\noindent Moreover,

$$
\log\delta_G=
34\alpha_1+66\alpha_2+96\alpha_3+52\alpha_4.
$$

Just as before, it follows that $\pi_{\min}$ 
is strongly $L^{p+\eps}$ if and only if 
$p\geq {26\over 3}$.

\end{enumerate}

\section{Minimal decay and isolatedness of $\pi_{\min}$}

In this section we complete the prove of Theorem 
\ref{mainthm}.
Let $G$ be one of the two groups introduced in 
the previous section.
Let $H_i, i\in\{1,2\}$ be defined as follows.\\

\noindent$\diamond \ H_1$ is the rank one subgroup corresponding to 
the highest root $\beta$.\\

\noindent $\diamond \ H_2$ is the subgroup corresponding to
all roots perpendicular to the highest root (i.e. its relative
diagram is the subdiagram $\mathbf C_3$ of $\mathbf F_4$).\\

\begin{lemma}
Let $G, H_i\ \  (i\in\{1,2\})$ be as above.
Let $\pi$ be an irreducible unitary representation
of $G$ which is not the
trivial or the minimal representation.
Then $\pi_{|H_1}$ is strongly $L^{2+\eps}$
and $\pi_{|H_2}$ is strongly $L^{4+\eps}$. 

\end{lemma}

\begin{proof}
For $i=1$, this follows from \cite[Theorem 4.1]{Oh}.
We prove the lemma for $i=2$. \\

As before, let $P=LN$ be
the Heisenberg parabolic of $G$.
Using Mackey
theory as 
in (\ref{mackeyPP}),
$\pi_{|[P,P]}$ can be expressed
as a direct integral of representations of the form
$$
\Res_{[P,P]}^P\Ind_{[P,P]}^P(\nu\otimes\rho).
$$
Let $M=[L,L]$. 
Since $\pi$ is not minimal,
it follows from Proposition 
\ref{UNIQUENESSOF}
that 
$\pi$ has to be of rank at least two.Therefore 
$\pi_{|M}$ is
a direct integral of representations of the form

\begin{equation}
\label{directintegralspectrum}
\nu_\chi\otimes\rho'_\chi
\end{equation}

\noindent where $\nu_\chi$ is a
representation of $M$ of rank at least one,
i.e. it does not have 
a nonzero $M$-fixed vector, and $\rho'_\chi$ is 
given as
$$
\rho'_\chi=\Res_{M}^{[P,P]} \rho_\chi
$$ 
where $\rho_\chi$ is the extension of the representation
of $N$ with central character $\chi$ to $M$ using the construction
of the oscillator representation \cite{Weil}. Although the 
central character $\chi$ can be arbitrary, there will be only
finitely many representations $\rho'_\chi$ (at most equal to
the number of elements of $F^\times/(F^\times)^2$). Therefore
$\pi_{|M}$ is a finite direct sum of representations
of the form given in (\ref{directintegralspectrum}),
and moreover since $\pi$ is of rank at least two, none of
the $\nu_\chi$'s can have trivial subrepresentations. \\

Suppose $\nu_\chi$ is strongly $L^{p+\eps}$ and $\rho'_\chi$ is strongly
$L^{q+\eps}$. 
Then 
by an application of H\"{o}lder's inequality,
one can see that $\nu_\chi\otimes\rho'_\chi$
is strongly $L^{r+\eps}$ where ${1\over r}={1\over p}+{1\over q}$, and
if this is true for all the direct summands of $\pi$ 
then $\pi_{|M}$ will be strongly $L^{r+\eps}$ as well. Therefore what remains
to be done is to find suitable values of $p$ and $q$ such that 
${1\over p}+{1\over q}\geq{1\over 4}$. To this end we use 
\cite[Prop. 2.2]{LiZhu} followed by \cite[Lemma 7.3]{Oh}.\\

Let $\Phi(a)$
denote the function introduced in \cite[Prop. 2.2]{LiZhu}; i.e. $\Phi$
is an upper bound  on the matix coefficients of the set of $K$-finite
vectors in the restriction of the oscillator representation to $M$.\\

Let $\alpha_1,...,\alpha_4$ be simple roots of the Dynkin diagram
of $G$ numbered such that the highest root is 
$$
\beta=2\alpha_1+3\alpha_2+4\alpha_3+2\alpha_4.
$$

\noindent Let $a\in A_{H_2}$ and let $y_i=|e^{\alpha_i}(a)|$. 
Since $|e^\beta(a)|=1$, we have $|e^{\alpha_1}(a)|=y_2^{-{3\over 2}}y_3^{-2}
y_4^{-1}$. Therefore a simple calculation shows that
$$
\Phi(a)\leq y_2^{-3r-5\over 4}y_3^{-r-2}y_4^{-r-3\over 2}
$$
where $r$ is the dimension of a short root space.\\

Note that by \cite[lemma 7.3]{Oh}, we have:

\begin{center}
$\Diamond$\ 
$\rho'_\chi$ is strongly
$L^{p+\eps}$ if and only if the exponent of any $y_i$ in
$\Phi^p(a)\delta_{H_2}(a)$ is nonpositive.\\
\end{center}

On the other hand,

$$\delta_{H_2}(a)=y_2^{3r+3}y_3^{6r+4}y_4^{4r+2}.$$

Therefore 

\begin{equation}
\notag
\delta_{H_2}(a)=\left\{
\begin{array}{lc}
y_2^9y_3^{16}y_4^{10} & \textrm{\ if}\ G \textrm{\ is}\ 
\!\!\ ^2\mathbf E_{6,4}^2 \\
 & \\
y_2^{15}y_3^{28}y_4^{18} & \textrm{if}\ G \textrm{\ is}\ 
\mathbf E_{7,4}^9 \\
\end{array}
\right.
\end{equation}

\noindent and consequently,

$$
q=\left\{
\begin{array}{lc}
4 & \textrm{\ if}\ G\ \textrm{\ is}\ \!\!\ ^2\mathbf E_{6,4}^2\\
{36\over 7} & \textrm{if}\ G\ \textrm{\ is}\ \mathbf E_{7,4}^9
\end{array}
\right.
$$

For $G=\!\!\ ^2\mathbf E_{6,4}^2$ the lemma is already proved,
since one can definitely take $p<\infty$ (because
$H_2$ has property $T$).
For $ \mathbf E_{7,4}^9$, 
it follows from
\cite[Theorem 7.4]{Oh} that $p\leq 18$. Note that 
although Oh's result is stated for irreducible unitary
representations, it is easily seen to hold for {\it any}
unitary representation without a nonzero 
fixed vector.
It is not hard to check that
${7\over 36}+{1\over 18}={1\over 4}$! This proves the lemma for
$i=2$.\\

\end{proof}

We will prove Theorem \ref{mainthm}
now.
We use Proposition \ref{asledecay}.\\

Let $a\in A_G$ and let $y_i=|e^{\alpha_i}(a)|$.\\

\begin{enumerate}
\item[$\bullet {\ ^2\bf E_{6,4}^2}: $]
We have\\
$$\delta_G(a)=y_1^{22}y_2^{42}y_3^{60}y_4^{32}
\textrm{\ and\ } 
\delta_{H_1}^{1\over 2}(a)
\delta_{H_2}^{1\over 4}(a)=
y_1y_2^{15\over 4}y_3^6y_4^{7\over 2}.$$\\

\noindent Consider the
element $w=s_{\alpha_4}s_{\alpha_2}s_
{\alpha_3}s_{\alpha_2}s_{\alpha_1}$
of the Weyl group of the relative root system
of $G$. 
Let $H_i^w$ ($i\in\{1,2\}$) 
denote the conjugate of $H_i$ by
$w$, i.e. its root system consists of all roots
$w\cdot \gamma$ where $\gamma$ is in the 
root system of $H_i$.
Then \\

$\delta_{H_1^w}^{1\over 2}(a)
\delta_{H_2^w}^{1\over 4}(a)=
y_1^{11\over 4}y_2^{21\over 4}y_3^{15\over 2}y_4^4$.\\

\noindent Proposition \ref{asledecay} implies that
any nontrivial and 
non-minimal irreducible representation of 
$G$ is strongly $L^{8+\eps}$; i.e.
$p(G)\leq 8$. By the result of section 6 we should
have $p(G)=8$.\\

\item[$\bullet {\bf E_{7,4}^9}: $] We have\\
$$\delta_G(a)=y_1^{34}y_2^{66}y_3^{96}y_4^{52}
\textrm{\ and\ } 
\delta_{H_1}^{1\over 2}\delta_{H_2}^{1\over 4}(a)=y_1y_2^{21\over 4}
y_3^{9}y_4^{11\over 2}.$$\\

\noindent Again we conjugate the groups by 
$w=s_{\alpha_4}s_{\alpha_1}s_{\alpha_2}s_{\alpha_3}
s_{\alpha_2}s_{\alpha_1}$ 
and we
get\\

$\delta_{H_1^w}^{1\over 2}(a)\delta_{H_2^w}^{1\over 4}(a)
=y_1^{9\over 2}y_2^{35\over 4}y_3^{25\over 2}y_4^7$\\

\noindent which gives $r\leq {192\over 25}<8$. Therefore, any 
non-minimal
irreducible representation of $G$ is strongly $L^{8+\eps}$,
and therefore $(\Chi^{1\over 4},\dim(\nu))$-bounded. However,
by Proposition  \ref{p6},
the minimal representation of $G$ can not be 
$(\Chi^{1\over 4},\dim(\nu))$-bounded since it is not strongly
$L^{8+\eps}$. Therefore the complement of the set
$\{{\mathrm{trivial}},\pi_{\min}\}$ is a closed set in the unitary dual
of $G$; i.e. $\pi_{\min}$ is isolated as well. Moreover, by section 6,
$p(G)={26\over 3}$.

\end{enumerate}

\end{document}